# $\bar{\partial}$-problem in fiber bundles for decreasing $(0,1)$-forms

Małgorzata Urlińska

**Abstract**

In this paper we consider the $\bar{\partial}$-problem in fiber bundles (fibers biholomorphic to $\mathbb{C}^k$, $k \geqslant 1$), namely the equation $\bar{\partial}\sigma = \omega$ for $(0,1)$-forms $\omega$ which decrease along the fibers. The order of decrease is slightly more than one. The important fact is that we do not assume that $\omega$ has compact support. The main theorem says that the equation has a solution which also decreases along fibers, however, not necessarily with the order as the original form. Existence of solution of the above mentioned $\bar{\partial}$-problem can be applied in various situations in Complex Analysis, in particular, to the Hartogs extension phenomenon.

## 1 Introduction

Inhomogeneous Cauchy-Riemann equation or $\bar{\partial}$-problem, namely the equation $\bar{\partial}\sigma = \omega$, where $\omega$ is a closed $(p,q)$ form, is a basic problem in Complex Analysis. Moreover, $\bar{\partial}$-problem has deep consequences on Algebraic Geometry, Partial Differential Equations and other areas of mathematics. The existence of solution of the $\bar{\partial}$-problem has applications to integral formulas, holomorphic extension, holomorphic approximations or vanishing of cohomology groups.

Usually for an explicit solution of the equation $\bar{\partial}\sigma = \omega$, say in $\mathbb{C}^n$, where $\omega$ is a $(0,1)$-closed form, we need that the support of $\omega$ is compact. If the support of $\omega$ is not compact, then we need additional assumptions on the domain where the equation is considered, for instance, Runge domains (see [14]). In Runge domains, we know that a solution of the equation $\bar{\partial}\sigma = \omega$ exists, but since the solution is not given explicitly, we cannot control some properties of $\sigma$ in terms of properties of $\omega$.

In this paper we will consider the equation $\bar{\partial}\sigma = \omega$ in fiber bundles $X$ and in the case of $(0,1)$-forms $\omega$ which decrease along the fibers. Roughly speaking, a fiber bundle is a complex manifold $X$ together with a holomorphic projection $\pi$ onto another complex manifold $\mathcal{B}$, called the base, and we will be assuming that the fibers $\pi^{-1}(p)$ are biholomorphic to $\mathbb{C}^k$ ($k = 1, 2, \ldots$ is fixed). Moreover, locally a chart $X_s$, $X_s \subset X$, is biholomorphic with $U_s \times \mathbb{C}^k$, where $U_s$ is a domain in $\mathbb{C}^n$, using the coordinates $(z, w) = (z_1, \ldots, z_n, w_1, \ldots, w_k)$. In the paper [7] a similar problem was considered but for forms $\omega(z, w)$ which have compact support along the fibers. In various situations, the compactness of $\text{supp}\,\omega$ is a too restrictive assumption. In this paper we drop the condition of compactness of support along the fibers and replace it by a much weaker condition that the form is decreasing along the fibers. Roughly speaking (more precisely it is defined in Section 2), in each chart $X_s$, the form

$$\omega(z,w) = a_1(z,w)\,d\bar{z}_1 + \ldots + a_n(z,w)\,d\bar{z}_n + b_1(z,w)\,d\bar{w}_1 + \ldots + b_k(z,w)\,d\bar{w}_k$$

is decreasing along the fibers if the coefficients of the form are of class $C^1$ and there exist $\varepsilon > 0$ and a constant $C > 0$ such that the coefficients satisfy the conditions

$$\lim_{|w|\to\infty} a_\alpha(z,w) = 0, \quad |b_\gamma(z,w)| \leqslant \frac{C}{1+\|w\|^{1+\varepsilon}}, \quad |Db_\gamma(z,w)| \leqslant \frac{C}{1+\|w\|^{1+\varepsilon}}, \tag{1}$$

where $D$ is any partial derivative of first order with respect to real variables $x_\beta = \operatorname{Re} z_\beta$ or $y_\beta = \operatorname{Im} z_\beta$ or $u_\delta = \operatorname{Re} w_\delta$ or $v_\delta = \operatorname{Im} w_\delta$, $\beta = 1,\ldots,n$, $\delta = 1,\ldots,k$. Here, $\|w\|^{1+\varepsilon}$ is a short notation for $|w_1|^{1+\varepsilon} + \ldots + |w_k|^{1+\varepsilon}$, but $|w_\delta|$ is the standard module of a complex number, $|w| = \sqrt{|w_1|^2 + \ldots + |w_k|^2}$. Of course, it is possible to choose other, almost equivalent conditions to (1), however, it seems that (1) is the most convenient and "symmetric" with respect to variables $w_1, \ldots, w_k$.

The main result can be stated (more precise version is in Theorem 2, beginning of Section 3):

**Theorem 1 (Main Theorem).** *Let $X$ be a fiber bundle and let $\omega(z,w)$ be a $(0,1)$-form that vanishes along the fibers. Then there exists a globally defined $C^1$-function $B(z,w)$ on $X$ such that $\lim_{|w|\to\infty} B(z,w) = 0$ (in each coordinate chart) and $\bar\partial B = \omega$.*

The theorems like Theorem 1 or from [7] are very useful in the Hartogs-type extension (see [13], [12]) of holomorphic or Cauchy-Riemann functions in a wide class of complex manifolds, see [15], [16], [18], [2], [3], [4], [5], [6], [9], [10]. Applications are given in [7] and we do not repeat them here. The advantage of the Main Theorem over the results in [7] is that it will allow us to apply the existence of extension to functions or forms on complex manifolds which are not necessarily holomorphic outside some set. It is enough that the functions or forms vanish at a point of the manifold or vanish on a complex submanifold, and after removing the point or the submanifold, we get a fiber-type bundle, like in [6] or toric manifolds [8]. As an illustration of this, we provide an example of such situation in the last Section 4.

## 2 Definitions and notation

**Fiber bundles.** Let $X$ be a complex connected manifold of complex dimension $n+k$, where $n, k = 1, 2, \ldots$. Assume that $X$ has a locally finite, open cover $\{X_s\}_{s \in S}$, where $X_s$ is biholomorphically equivalent to $U_s \times \mathbb{C}^k$, where $U_s$ is a domain in $\mathbb{C}^n$. The coordinates on $X_s$ are of the form

$$(z, w) = (z_s, w_s) = (z_{s1}, \ldots, z_{sn}, w_{s1}, \ldots, w_{sk}) = (z_1, \ldots, z_n, w_1, \ldots, w_k)$$

and they give a biholomorphism from $X_s$ onto $U_s \times \mathbb{C}^k$. The real and imaginary parts of the coordinates are

$$z_{s\alpha} = x_{s\alpha} + i y_{s\alpha} = x_\alpha + i y_\alpha, \quad \alpha = 1,\ldots,n, \qquad w_{s\gamma} = u_{s\gamma} + i v_{s\gamma} = u_\gamma + i v_\gamma, \quad \gamma = 1,\ldots,k.$$

The transition functions are of the form
$$(z_s, w_s) = (f_{st}(z_t), g_{st}(z_t, w_t)), \quad s, t \in S,$$
the function $f_{st}$ depends on $z_t$ only, and moreover the second mapping $g_{st}$ gives a biholomorphism
$$\mathbb{C}^k \ni w_t \longrightarrow g_{st}(z_t, w_t) \in \mathbb{C}^k.$$

We also use the following notation,
$$w = (w_1, \ldots, \ldots, w_k)$$
$$|w| = \left(|w_1|^2 + \ldots + |w_k|^2\right)^{1/2} = \sqrt{|w_1|^2 + \ldots + |w_k|^2}$$
$$w' = (w_1, \ldots, w_{\delta-1}, w_{\delta+1}, \ldots, w_k)$$
$$\|w'\|^{1+\varepsilon} = |w_1|^{1+\varepsilon} + \ldots + |w_{\delta-1}|^{1+\varepsilon} + |w_{\delta+1}|^{1+\varepsilon} + \ldots + |w_k|^{1+\varepsilon}$$
$$\|w'\| = |w_1| + \ldots + |w_{\delta-1}| + |w_{\delta+1}| + \ldots + |w_k|$$
Purposely we do not put the index $\delta$ at $w'$, hopefully there will be no confusion.

**(0,1)-forms.** In this paper we will deal only with functions and $(0,1)$-forms. The functions and the forms will be defined on the fiber bundle $X$ or in the local coordinate domain $X_s$. A $(0,1)$-form defined on $X_s$ can be written

$$\omega_s(z_s, w_s) = \omega(z, w) = \sum_{\alpha=1}^{n} a_\alpha(z, w)\, d\bar{z}_\alpha + \sum_{\gamma=1}^{k} b_\gamma(z, w)\, d\bar{w}_\gamma. \tag{2}$$

Always we will be assuming that the form is of class $C^1$, i.e., the coefficients $a_\alpha(z, w)$ and $b_\gamma(z, w)$ are $C^1$ in the domain when the form is considered. We have a natural split of the form into two parts, the so-called $a$-part (or $z$-part) and $b$-part (or $w$-part):

$$a(z, w) = \sum_{\alpha=1}^{n} a_\alpha(z, w)\, d\bar{z}_\alpha, \qquad b(z, w) = \sum_{\gamma=1}^{k} b_\gamma(z, w)\, d\bar{w}_\gamma. \tag{3}$$

The form $\omega$ is $\bar{\partial}$-closed if $\bar{\partial}\omega = 0$, which is equivalent to the compatibility conditions

$$\frac{\partial a_\alpha}{\partial \bar{z}_\beta} = \frac{\partial a_\beta}{\partial \bar{z}_\alpha}, \quad \frac{\partial a_\alpha}{\partial \bar{w}_\gamma} = \frac{\partial b_\gamma}{\partial \bar{z}_\alpha}, \quad \frac{\partial b_\gamma}{\partial \bar{w}_\delta} = \frac{\partial b_\delta}{\partial \bar{w}_\gamma} \quad \text{where} \quad \alpha, \beta = 1, \ldots, n, \quad \gamma, \delta = 1, \ldots, k. \tag{4}$$

**Decreasing forms.** Let $\omega$ be a $(0,1)$-form defined on a fiber bundle $X$. We say that the form is decreasing along the fibers if in each coordinate domain $X_s \approx U_s \times \mathbb{C}^k$, for any compact set $K \subset U_s$, the coefficients satisfy the conditions:

*There exist $\varepsilon = \varepsilon(s, K) > 0$ and a constant $C = C(s, K, \varepsilon)$ such that*

$$\lim_{|w| \to \infty} a_\alpha(z, w) = 0 \qquad \text{uniformly with respect to } z \in K, \quad \alpha = 1, \ldots, n;$$

$$|b_\gamma(z, w)| \leq \frac{C}{1 + |w_1|^{1+\varepsilon} + \ldots + |w_k|^{1+\varepsilon}} \qquad \text{uniformly with respect to } z \in K, \quad \gamma = 1, \ldots, k; \qquad (5)$$

$$|Db_\gamma(z, w)| \leq \frac{C}{1 + |w_1|^{1+\varepsilon} + \ldots + |w_k|^{1+\varepsilon}} \qquad \text{uniformly with respect to } z \in K, \quad \gamma = 1, \ldots, k;$$

where $D$ is any partial derivative

$$\frac{\partial}{\partial x_\beta}, \quad \frac{\partial}{\partial y_\beta}, \quad \frac{\partial}{\partial u_\delta}, \quad \frac{\partial}{\partial v_\delta}, \qquad \beta = 1, \ldots, n, \quad \delta = 1, \ldots, k.$$

Here the uniform on $K$ limit means, e.g. for the first limit above,

$$\forall \vartheta > 0 \ \exists M_{K,\vartheta} > 0 \ \forall |w| \geq M_{K,\vartheta} \ \forall z \in K \quad |a_\alpha(z, w)| < \vartheta.$$

The number $\varepsilon > 0$ in (5) can be arbitrarily small. When $\varepsilon$ is chosen smaller, then the conditions in (5) become weaker.

We note that if a closed $(0, 1)$-form $\omega(z, w)$ is decreasing along the fibers, then the $a$-part of the form (as in (3)) is uniquely determined by the $b$-part. Namely we have: suppose that $b(z, w) \equiv 0$, then from the compatibility conditions (4) we obtain

$$\frac{\partial a_\alpha}{\partial \bar{w}_\gamma}(z, w) \equiv 0 \qquad \alpha = 1, \ldots, n, \quad \gamma = 1, \ldots, k,$$

which gives that the coefficients $a_\alpha(z, w)$ are holomorphic with respect to $w$. Since these coefficients are decreasing to zero with respect to $w$ as $|w| \to \infty$, we have that $a_\alpha(z, w) \equiv 0$.

## 3 Solution of the $\bar{\partial}$-problem and proof of the theorem

The main goal of this section is to prove the following

**Theorem 2** *Let $X$ be a fiber bundle and let $\omega = \omega(z, w)$ be a $(0, 1)$-form that vanishes along the fibers, as it is defined in (5). Then there exists a globally defined $C^1$-function $B(z, w)$ on $X$ such that*

1. *$\bar{\partial} B = \omega$.*

2. In each coordinate map $X_s$, as defined in Section 2 (other notation also from this section), the solution satisfies the estimate

$$|B(z,w)| \leq \text{const} \int_{\mathbb{C}} \frac{1}{|\zeta|} \cdot \frac{1}{\left(1 + \|w'\|^{1+\varepsilon} + \big||w_\delta| + \zeta\big|^{1+\varepsilon}\right)} |d\zeta \wedge d\bar\zeta| \quad \text{for any} \quad \delta = 1, \ldots, k,$$

when $z$ is running over a compact set and the constant depends on this compact set.

3. The function $B(z,w)$ is bounded in each coordinate chart, uniformly when $z$ is running over a compact set.

4. $\lim_{|w_\delta| \to \infty} B(z,w) = 0$, $\delta = 1, \ldots, k$, in each coordinate map and the limit is uniform when $z$ is running over a compact set, and $w' = (w_1, \ldots, w_{\delta-1}, w_{\delta+1}, \ldots, w_k)$ is any.

## 3.1 Definition of the function $B_\delta(z,w)$

Our goal is to solve the equation $\bar\partial \sigma = \omega$ in the fiber bundle $X$. We have to be careful with the method (as in [7] of solving the problem since the functions $b_\delta(z,w)$, $\delta = 1, \ldots, k$, from (2) do not have compact support with respect to $w$. We define

$$B_\delta(z,w) = \frac{1}{2\pi i} \int_{\mathbb{C}} \frac{b_\delta(z, w_1, \ldots, w_{\delta-1}, \zeta, w_{\delta+1}, \ldots, w_k)}{\zeta - w_\delta} d\zeta \wedge d\bar\zeta \tag{6}$$

$$= \frac{1}{2\pi i} \int_{\mathbb{C}} \frac{b_\delta(z, w_1, \ldots, w_{\delta-1}, w_\delta + \zeta, w_{\delta+1}, \ldots, w_k)}{\zeta} d\zeta \wedge d\bar\zeta \tag{7}$$

If $\zeta = \xi + i\eta$, then

$$d\zeta \wedge d\bar\zeta = (d\xi + id\eta) \wedge (d\xi - id\eta) = d\xi \wedge d\xi - id\xi \wedge d\eta + id\eta \wedge d\xi + d\eta \wedge d\eta = -2id\xi \wedge d\eta.$$

Later on, $|d\zeta \wedge d\bar\zeta|$ will be understood as

$$|d\zeta \wedge d\bar\zeta| = |-2id\xi \wedge d\eta| = 2\, d\xi\, d\eta.$$

We note that the integrals (6) or (7) make sense because the function $b_\delta$ is decreasing with respect to the variables $w_1, \ldots, w_{\delta-1}, \zeta, w_{\delta+1}, \ldots, w_k$. Namely we have:

Suppose that $w_\delta$ is running over a compact set $K_\delta \subset \mathbb{C}$. Then from the assumption on the decreasing condition, there exists $\varepsilon > 0$ and a constant $C$ such that

$$\left| \frac{b_\delta(z, w_1, \ldots, w_{\delta-1}, w_\delta + \zeta, w_{\delta+1}, \ldots, w_k)}{\zeta} \right| \leq$$

$$\leq \frac{C}{|\zeta| \cdot (1 + |w_1|^{1+\varepsilon} + \ldots + |w_{\delta-1}|^{1+\varepsilon} + |w_\delta + \zeta|^{1+\varepsilon} + |w_{\delta+1}|^{1+\varepsilon} + \ldots + |w_k|^{1+\varepsilon})} \leq$$

$$\leq \frac{C}{|\zeta| \cdot (1 + |w_\delta + \zeta|^{1+\varepsilon})}$$

Since $w_\delta$ is running over a compact set $K_\delta \subset \mathbb{C}$, possibly increasing the constant $C$ to $\widetilde{C}$, the last quotient can be estimated by

$$\frac{\widetilde{C}}{|\zeta| \cdot (1 + |\zeta|^{1+\varepsilon})}$$

so we have

$$\left| \frac{b_\delta(z, w_1, \ldots, w_{\delta-1}, w_\delta + \zeta, w_{\delta+1}, \ldots, w_k)}{\zeta} \right| \leq \frac{\widetilde{C}}{|\zeta| \cdot (1 + |\zeta|^{1+\varepsilon})}. \tag{8}$$

Clearly the integral of the function on the right above exists and is finite. Namely we have

$$\left| \int_\mathbb{C} \frac{1}{|\zeta| \cdot (1 + |\zeta|^{1+\varepsilon})} d\zeta \wedge d\bar\zeta \right| \leq \int_\mathbb{C} \frac{1}{|\zeta| \cdot (1 + |\zeta|^{1+\varepsilon})} |d\zeta \wedge d\bar\zeta|$$

$$\stackrel{\zeta = re^{i\theta}}{=} \int_0^{2\pi} \int_0^\infty \frac{1}{r(1 + r^{1+\varepsilon})} 2r \, dr \, d\theta$$

$$= 4\pi \int_0^\infty \frac{1}{1 + r^{1+\varepsilon}} dr$$

$$= 4\pi \int_0^1 \frac{1}{1 + r^{1+\varepsilon}} dr + 4\pi \int_1^\infty \frac{1}{1 + r^{1+\varepsilon}} dr$$

$$\leq 4\pi \int_0^1 dr + 4\pi \int_1^\infty \frac{1}{r^{1+\varepsilon}} dr$$

$$= 4\pi + 4\pi \left[ -\frac{1}{\varepsilon} \cdot \frac{1}{r^\varepsilon} \right]_1^\infty$$

$$= 4\pi + 4\pi \cdot \frac{1}{\varepsilon}$$

$$= 4\pi \left( 1 + \frac{1}{\varepsilon} \right)$$

## 3.2 Estimates of the function $B_\delta(z, w)$ along the fibers

We investigate the behavior of the function $B_\delta(z, w)$ along the fibers. By assumption, the function $b_\delta(z, w)$ and its first order derivatives satisfy the conditions (5). We will estimate the function $B_\delta(z, w)$ defined in (6) or (7). We use the notation as in the beginning of Section 2.

### 3.2.1 First estimate of $B_\delta(z, w)$

In this subsection we will prove the following

**Lemma 1** *With all the notation as above, we have the estimate*

$$|B_\delta(z,w)| \leqslant \int_{\mathbb{C}} \frac{1}{|\zeta|} \cdot \frac{C}{(1 + \|w'\|^{1+\varepsilon} + |w_\delta + \zeta|^{1+\varepsilon})} \, |d\zeta \wedge d\bar\zeta|$$

$$\leqslant \int_{\mathbb{C}} \frac{1}{|\zeta|} \cdot \frac{C}{(1 + \|w'\|^{1+\varepsilon} + |\zeta|^{1+\varepsilon})} \, |d\zeta \wedge d\bar\zeta|$$

$$\leqslant \int_{\mathbb{C}} \frac{1}{|\zeta|} \cdot \frac{C}{(1 + |\zeta|^{1+\varepsilon})} \, |d\zeta \wedge d\bar\zeta|$$

*where the constant $C$ is the same as in the assumption (5) of the decreasing condition of coefficients of the $(0,1)$-form $b(z, w)$.*

**Remark 1** *We note that the second inequality above holds between the integrals, not as implication of a corresponding inequality between the integrands. Actually the inequality between the integrands is not true.*

**Proof Lemma 1.** We have the following estimates:

$$|B_\delta(z,w)| = \left| \frac{1}{2\pi i} \int_{\mathbb{C}} \frac{b_\delta(z, w_1, \ldots, w_{\delta-1}, w_\delta + \zeta, w_{\delta+1}, \ldots, w_k)}{\zeta} d\zeta \wedge d\bar\zeta \right|$$

$$\leqslant \frac{1}{2\pi} \int_{\mathbb{C}} \frac{|b_\delta(z, w_1, \ldots, w_{\delta-1}, w_\delta + \zeta, w_{\delta+1}, \ldots, w_k)|}{|\zeta|} \, |d\zeta \wedge d\bar\zeta|$$

$$\leqslant \frac{1}{2\pi} \int_{\mathbb{C}} \frac{1}{|\zeta|} \cdot \frac{C}{(1 + |w_1|^{1+\varepsilon} + \ldots + |w_{\delta-1}|^{1+\varepsilon} + |w_\delta + \zeta|^{1+\varepsilon} + |w_{\delta+1}|^{1+\varepsilon} + \ldots + |w_k|^{1+\varepsilon})} \, |d\zeta \wedge d\bar\zeta|$$

$$= \frac{1}{2\pi} \int_{\mathbb{C}} \frac{1}{|\zeta|} \cdot \frac{C}{(1 + \|w'\|^{1+\varepsilon} + |w_\delta + \zeta|^{1+\varepsilon})} \, |d\zeta \wedge d\bar\zeta|$$

In the next sequence of estimates we switch to polar coordinates,

$$2\pi |B_\delta(z,w)| \leq \int_{\mathbb{C}} \frac{1}{|\zeta|} \cdot \frac{C}{(1+\|w'\|^{1+\varepsilon}+|w_\delta+\zeta|^{1+\varepsilon})} |d\zeta \wedge d\bar\zeta| \qquad (9)$$

$$\stackrel{\zeta=re^{i\theta}}{=} \int_0^{2\pi} \int_0^\infty \frac{1}{|re^{i\theta}|} \cdot \frac{C}{(1+\|w'\|^{1+\varepsilon}+|w_\delta+re^{i\theta}|^{1+\varepsilon})} 2r\, dr\, d\theta$$

$$= \int_0^{2\pi} \int_0^\infty \frac{C}{(1+\|w'\|^{1+\varepsilon}+|w_\delta+re^{i\theta}|^{1+\varepsilon})} 2dr\, d\theta$$

$$= \int_0^{2\pi} \int_0^\infty \frac{C}{1+\|w'\|^{1+\varepsilon}+\big||w_\delta|+re^{i\theta}\big|^{1+\varepsilon}} 2dr\, d\theta$$

Obviously the last integral depends on $\|w'\|$; also it depends on $|w_\delta|$ only because of the integration with respect to $\theta$. Moreover it is clear that in the case $k \geq 2$ if $w_\delta$ is running over a compact set, then the integral is tending to zero when $\|w'\| \to \infty$; if $k = 1$, the term $\|w'\|^{1+\varepsilon}$ is absent.

Now we look at the integral (9) (after dropping the constant $C$):

$$\int_{\mathbb{C}} \frac{1}{|\zeta|} \cdot \frac{1}{(1+\|w'\|^{1+\varepsilon}+|w_\delta+\zeta|^{1+\varepsilon})} |d\zeta \wedge d\bar\zeta| \qquad (10)$$

We ask the question:

> With fixed $w'$, for what $w_\delta$ the above integral takes the biggest value?

If we look at this integral geometrically, the answer is:

> For such $w_\delta$ at which $|w_\delta + \zeta|$ is the smallest when $\zeta$ is zero, that is, $w_\delta = 0$. Another explanation of this property is that we want the fraction $\dfrac{1}{|\zeta|}$ to take the maximum value while at the same time $|w_\delta + \zeta|$ takes the smallest value.

If the geometric argument is correct, we have the inequality

$$\int_{\mathbb{C}} \frac{1}{|\zeta|} \cdot \frac{1}{(1+\|w'\|^{1+\varepsilon}+|w_\delta+\zeta|^{1+\varepsilon})} |d\zeta \wedge d\bar\zeta| \leq \int_{\mathbb{C}} \frac{1}{|\zeta|} \cdot \frac{1}{(1+\|w'\|^{1+\varepsilon}+|\zeta|^{1+\varepsilon})} |d\zeta \wedge d\bar\zeta| < +\infty \qquad (11)$$

Because the geometric explanation is heuristic, now we look at the above inequality more algebraically. Since we know that the integral (10) depends only on $|w_\delta|$, let us examine the function

$$F(x) = \int_{\mathbb{C}} \frac{1}{|\zeta|} \cdot \frac{1}{(1+\|w'\|^{1+\varepsilon}+|x+\zeta|^{1+\varepsilon})} |d\zeta \wedge d\bar\zeta| \qquad \text{for } x \geq 0 \qquad (w' \text{ is fixed}). \qquad (12)$$

Our near goal is to calculate $F'(x)$, $x > 0$, however, at this moment we do not know if we can enter with the derivative under the integral sign. First we calculate the derivative of the integrand (actually of the second factor) and estimate it. Since the derivative is a local notion, we assume that $x$ is running over a compact interval $I$, $I \subset (0, \infty)$. We have

$$\left| \frac{\partial}{\partial x} \left[ \frac{1}{(1 + \|w'\|^{1+\varepsilon} + |x + \zeta|^{1+\varepsilon})} \right] \right| = \frac{1}{(1 + \|w'\|^{1+\varepsilon} + |x + \zeta|^{1+\varepsilon})^2} \cdot \frac{(1+\varepsilon)|x+\xi|}{[(x+\xi)^2 + \eta^2]^{\frac{1}{2} - \frac{\varepsilon}{2}}}$$

$$\leqslant \frac{1}{(1 + \|w'\|^{1+\varepsilon} + |x + \zeta|^{1+\varepsilon})^2} \cdot \frac{(1+\varepsilon)|x+\xi|}{[(x+\xi)^2]^{\frac{1}{2}}} \cdot \left[(x+\xi)^2 + \eta^2\right]^{\frac{\varepsilon}{2}}$$

$$= \frac{1+\varepsilon}{(1 + \|w'\|^{1+\varepsilon} + |x + \zeta|^{1+\varepsilon})^2} \cdot \left[(x+\xi)^2 + \eta^2\right]^{\frac{\varepsilon}{2}}$$

$$= \frac{1+\varepsilon}{(1 + \|w'\|^{1+\varepsilon} + |x + \zeta|^{1+\varepsilon})^2} \cdot |x+\zeta|^{\varepsilon}$$

$$\leqslant \frac{1+\varepsilon}{(1 + \|w'\|^{1+\varepsilon} + |x + \zeta|^{1+\varepsilon})}$$

$$\overset{x \in I}{\leqslant} \frac{(1+\varepsilon)\, C_I}{(1 + \|w'\|^{1+\varepsilon} + |\zeta|^{1+\varepsilon})}$$

where $C_I$ is a constant that depends on the interval $I$ only.

Consequently, we have an estimate of the derivative with respect to $x$ of the integrand in $F(x)$:

$$\frac{1}{|\zeta|} \left| \frac{\partial}{\partial x} \left[ \frac{1}{(1 + \|w'\|^{1+\varepsilon} + |x + \zeta|^{1+\varepsilon})} \right] \right| \leqslant \frac{1}{|\zeta|} \cdot \frac{(1+\varepsilon)\, C_I}{(1 + \|w'\|^{1+\varepsilon} + |\zeta|^{1+\varepsilon})} \qquad when \quad x \in I, \qquad (13)$$

for $\zeta \in \mathbb{C}$ and any fixed $w'$. This estimate allows us to enter with derivative of $F(x)$ under the integral sign. To be more precise, we can apply the Lebesgue Dominated Limit Theorem, see R. Sikorski [22] or R.G. Bartle [1] or H.L. Royden - P.M. Fitzpatrick [20], and a modified mean value theorem for vector-valued functions, see e.g., W. Rudin [21] p. 113, or S.M. Nikolsky [19] p. 119. We show all these applications in detail in Subsection 3.3, because of that, we do not repeat them here - the arguments are exactly the same.

We have,

$$F'(x) = -\int_{\mathbb{C}} \frac{1}{|\zeta|} \cdot \frac{\partial}{\partial x} \left[ \frac{1}{(1 + \|w'\|^{1+\varepsilon} + |x + \zeta|^{1+\varepsilon})} \right] |d\zeta \wedge d\bar{\zeta}| \qquad for \quad x > 0,$$

or equivalently,

$$
\begin{aligned}
F'(x) &= -\int_{\mathbb{C}} \frac{1}{|\zeta|} \cdot \frac{1}{(1+\|w'\|^{1+\varepsilon}+|x+\zeta|^{1+\varepsilon})^2} \cdot \frac{(1+\varepsilon)(x+\xi)}{[(x+\xi)^2+\eta^2]^{\frac{1}{2}-\frac{\varepsilon}{2}}} |d\zeta \wedge d\bar{\zeta}| \\
&= -\int_{\mathbb{C}} \frac{1}{[\xi^2+\eta^2]^{1/2}} \cdot \frac{1}{\left[1+\|w'\|^{1+\varepsilon}+[(x+\xi)^2+\eta^2]^{(1+\varepsilon)/2}\right]^2} \cdot \frac{2(1+\varepsilon)(x+\xi)}{[(x+\xi)^2+\eta^2]^{\frac{1}{2}-\frac{\varepsilon}{2}}} d\xi\, d\eta
\end{aligned}
$$

To simplify notation, we denote the integrand of the last integral by $G(\xi, \eta, x)$ (without the minus sign) and examine the behavior of $G$ with respect to $\xi$ with $\eta$, $x$ fixed. From Fig. 1, immediately we see

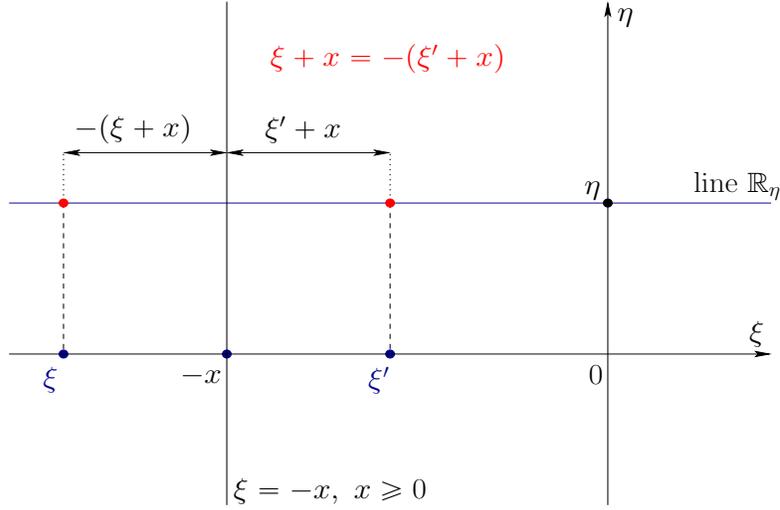

Figure 1: *Position of points at which to have an estimate*

that $\xi$ and $\xi'$ are positioned symmetrically with respect to the line passing through $(-x)$,

$$\xi' + x = |\xi + x|, \qquad \xi + x \leqslant 0, \qquad \xi' + x \geqslant 0, \qquad \xi^2 \geqslant (\xi')^2 \qquad (\text{recall that } x \geqslant 0),$$

and from these inequalities,

$$G(\xi, \eta, x) \leqslant 0, \qquad G(\xi', \eta, x) \geqslant 0, \qquad |G(\xi, \eta, x)| \leqslant G(\xi', \eta, x).$$

Using the above inequalities, we get that

$$-\int_{\mathbb{R}_\eta} G(\xi, \eta, x)\, d\xi \leqslant 0,$$

which implies that

$$F'(x) = -\int_{\mathbb{C}} G(\xi, \eta, x)\, d\xi\, d\eta \leqslant 0.$$

The inequality (11) is proved.

Using the same type of estimates as for the derivative $F'(x)$, we get that the function $F(x)$ is continuous on the interval $[0, \infty)$.

From (9) and (11), immediately we obtain

$$B_\delta(z, w) \quad \text{is bounded for all } z, w, \text{ and} \quad B_\delta(z, w) \longrightarrow 0 \quad \text{as} \quad \|w'\| \to \infty \quad (\text{if it is not void}), \tag{14}$$

however, we do not know its behavior with respect to $w_\delta$. So we will take a look at this case.

### 3.2.2 Second estimate of $B_\delta(z, w)$

Actually this subsection gives a geometric explanation of the property that the function $F(x)$ defined in (12) is decreasing. An additional property which we prove here is $\lim\limits_{|w_\delta| \to \infty} F(|w_\delta|) = 0$.

The main result in this section is the following

**Lemma 2** *The function $|B_\delta(z, w)|$ is bounded with respect to $w$ and uniformly for $z$ running in a compact set. Moreover, we have*

$$\begin{aligned}
|B_\delta(z, w)| &\leq C \int_\mathbb{C} \frac{1}{|\zeta|} \cdot \frac{1}{(1 + \|w'\|^{1+\varepsilon} + |w_\delta + \zeta|^{1+\varepsilon})} |d\zeta \wedge d\bar{\zeta}| \\
&= C \int_\mathbb{C} \frac{1}{|\zeta|} \cdot \frac{1}{\left(1 + \|w'\|^{1+\varepsilon} + \big||w_\delta| + \zeta\big|^{1+\varepsilon}\right)} |d\zeta \wedge d\bar{\zeta}| \quad (15) \\
&\longrightarrow 0 \quad \text{as} \quad |w_\delta| \to \infty,
\end{aligned}$$

*where the constant $C$ is the same as in the assumption (5) of the decreasing condition of coefficients of the $(0,1)$-form $b(z, w)$.*

**Proof Lemma 2.** We consider again the following integral

$$F(w_\delta) = \int_\mathbb{C} \frac{1}{|\zeta - w_\delta|} \cdot \frac{1}{(1 + \|w'\|^{1+\varepsilon} + |\zeta|^{1+\varepsilon})} |d\zeta \wedge d\bar{\zeta}| = \int_\mathbb{C} \frac{1}{|\zeta|} \cdot \frac{1}{(1 + \|w'\|^{1+\varepsilon} + |w_\delta + \zeta|^{1+\varepsilon})} |d\zeta \wedge d\bar{\zeta}| \tag{16}$$

From the considerations in the previous subsection, we know that $F$ depends on $|w_\delta|$ only, i.e., we have $F(w_\delta) = F(|w_\delta|)$. Geometrically it should be clear that $F(|w_\delta|) \to 0$ as $|w_\delta| \to \infty$, as is shown on Fig. 2. Namely, if we multiply the functions $\tau = \dfrac{1}{1 + \|w'\|^{1+\varepsilon} + |w_\delta + \zeta|^{1+\varepsilon}}$ and $\tau = \dfrac{1}{|\zeta|}$ and if $|w_\delta| \to \infty$, then it is obvious that the product becomes smaller and the "stack" gets thinner.

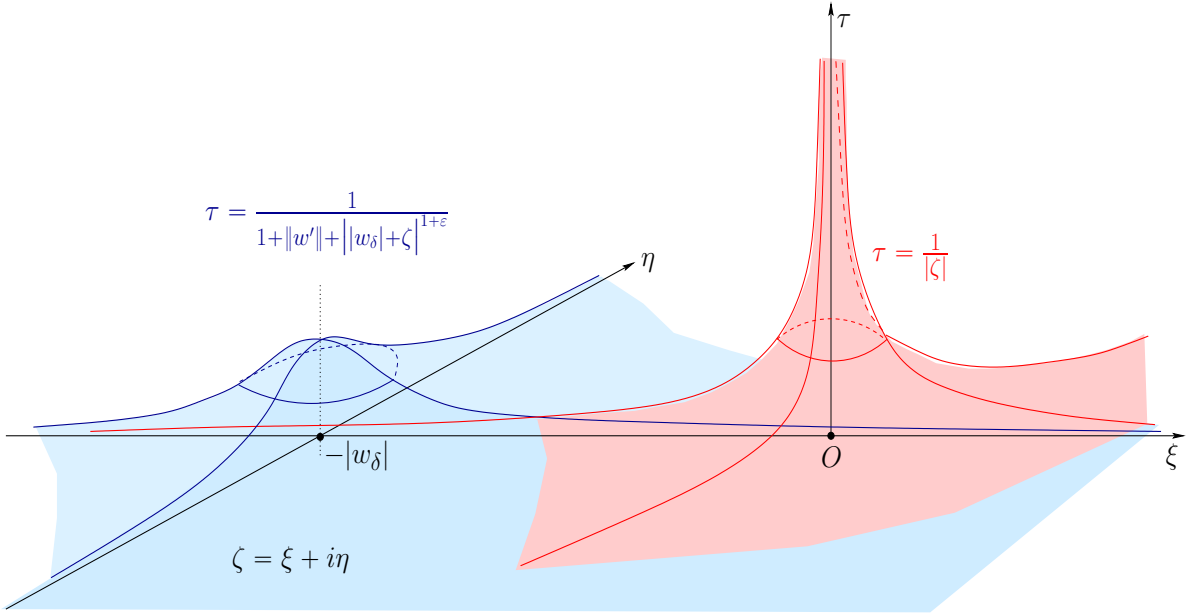

Figure 2: *Graphs of two functions*

Of course these arguments are not very precise, because of that we consider this integral (16) analytically. We have

$$\int_{\mathbb{C}} \frac{1}{|\zeta|} \cdot \frac{1}{(1+\|w'\|^{1+\varepsilon}+|w_\delta+\zeta|^{1+\varepsilon})} |d\zeta \wedge d\bar\zeta| = \int_0^{2\pi}\int_0^\infty \frac{1}{r} \cdot \frac{1}{(1+\|w'\|^{1+\varepsilon}+|w_\delta+re^{i\theta}|^{1+\varepsilon})} 2r\, dr\, d\theta$$

$$= \int_0^{2\pi}\int_0^\infty \frac{2}{1+\|w'\|^{1+\varepsilon}+|w_\delta+re^{i\theta}|^{1+\varepsilon}}\, dr\, d\theta$$

$$= \int_0^{2\pi}\int_0^\infty \frac{2}{1+\|w'\|^{1+\varepsilon}+\left||w_\delta|+re^{i\theta}\right|^{1+\varepsilon}}\, dr\, d\theta$$

For simplicity, denote $\rho = |w_\delta|$ and consider the function

$$g(\rho,\theta) = \int_0^\infty \frac{1}{1+\|w'\|^{1+\varepsilon}+|\rho+re^{i\theta}|^{1+\varepsilon}}\, dr$$

$$= \int_0^\infty \frac{1}{1+\|w'\|^{1+\varepsilon}+[(\rho+r\cos\theta)^2+(r\sin\theta)^2]^{(1+\varepsilon)/2}}\, dr$$

$$= \int_0^\infty \frac{1}{1+\|w'\|^{1+\varepsilon}+[\rho^2+2\rho r\cos\theta+r^2]^{(1+\varepsilon)/2}}\, dr$$

We will estimate the above integral. We have a sequence of equalities

$$g(\rho,\theta) = \int_0^\infty \frac{1}{1+\|w'\|^{1+\varepsilon}+[\rho^2+2\rho r\cos\theta+r^2]^{(1+\varepsilon)/2}}\,dr$$

$$= \int_0^\infty \frac{1}{1+\|w'\|^{1+\varepsilon}+[(r^2+2\rho r\cos\theta+\rho^2\cos^2\theta)-\rho^2\cos^2\theta+\rho^2]^{(1+\varepsilon)/2}}\,dr$$

$$= \int_0^\infty \frac{1}{1+\|w'\|^{1+\varepsilon}+[(r+\rho\cos\theta)^2+\rho^2(1-\cos^2\theta)]^{(1+\varepsilon)/2}}\,dr$$

$$= \int_0^\infty \frac{1}{1+\|w'\|^{1+\varepsilon}+[(r+\rho\cos\theta)^2+\rho^2\sin^2\theta]^{(1+\varepsilon)/2}}\,dr$$

$$\stackrel{r+\rho\cos\theta\leadsto r}{=} \int_{\rho\cos\theta}^\infty \frac{1}{1+\|w'\|^{1+\varepsilon}+[r^2+\rho^2\sin^2\theta]^{(1+\varepsilon)/2}}\,dr$$

Finally we estimate $B_\delta(z,w)$. Take any small positive number $\kappa$. From the interval $[0,2\pi]$ we remove three subintervals: $[0,\kappa]$, $[\pi-\kappa,\pi+\kappa]$, and $[2\pi-\kappa,2\pi]$. We denote the union of these intervals by $I_\kappa$. We can find a constant $M_\kappa>0$ such that

$$0\leqslant g(\rho,\theta)\leqslant \int_{-\infty}^\infty \frac{1}{1+\|w'\|^{1+\varepsilon}+[r^2+\rho^2\sin^2\theta]^{(1+\varepsilon)/2}}\,dr \leqslant \kappa \quad \text{for any}\quad \rho>M_\kappa,\ \theta\in[0,2\pi]\setminus I_\kappa. \tag{17}$$

If $\theta$ is in $I_\kappa$, then we have

$$0\leqslant g(\rho,\theta)\leqslant \int_{-\infty}^\infty \frac{1}{1+\|w'\|^{1+\varepsilon}+|r|^{1+\varepsilon}}\,dr \leqslant 2+\frac{2}{\varepsilon} \tag{18}$$

From (17) and (18) and using (16), we obtain

$$|B_\delta(z,w)|\leqslant C|F(w)| \leqslant C\int_{I_\kappa} g(\rho,\theta)\,d\theta + C\int_{[0,2\pi]\setminus I_\kappa} g(\rho,\theta)\,d\theta$$

$$\leqslant C\left[4\kappa\left(2+\frac{2}{\varepsilon}\right)+\kappa(2\pi-4\kappa)\right] \quad \text{if}\quad \rho=|w|>M_\kappa$$

We recall that $\varepsilon$ is a *fixed* positive number, therefore

$$C\kappa\left[4\left(2+\frac{2}{\varepsilon}\right)+(2\pi-4\kappa)\right]$$

can be made arbitrarily small if $\kappa$ is small.

## 3.3 Calculation of the derivative $\dfrac{\partial B_\delta}{\partial \bar{w}_\gamma}$

In this subsection we will prove the property that we can enter with differentiation under the integral sign. Since differentiation is a local notion, without any loss of generality we can assume that $w_1, \ldots, w_k$ are running over a compact set in $\mathbb{C}^k$.

**Lemma 3** *We have the property*

$$\begin{aligned}
\frac{\partial}{\partial \bar{w}_\gamma} B_\delta(z, w) &= \frac{1}{2\pi i} \frac{\partial}{\partial \bar{w}_\gamma} \int_\mathbb{C} \frac{b_\delta(z, w_1, \ldots, w_{\delta-1}, \zeta, w_{\delta+1}, \ldots, w_k)}{\zeta - w_\delta} d\zeta \wedge d\bar\zeta \quad (19) \\
&= \frac{1}{2\pi i} \frac{\partial}{\partial \bar{w}_\gamma} \int_\mathbb{C} \frac{b_\delta(z, w_1, \ldots, w_{\delta-1}, w_\delta + \zeta, w_{\delta+1}, \ldots, w_k)}{\zeta} d\zeta \wedge d\bar\zeta \\
&= \frac{1}{2\pi i} \int_\mathbb{C} \frac{1}{\zeta} \cdot \frac{\partial b_\delta}{\partial \bar{w}_\gamma}(z, w_1, \ldots, w_{\delta-1}, w_\delta + \zeta, w_{\delta+1}, \ldots, w_k) \, d\zeta \wedge d\bar\zeta
\end{aligned}$$

*The same is true for derivatives $\dfrac{\partial B_\delta}{\partial w_\gamma}(z, w)$, $\dfrac{\partial B_\delta}{\partial z_\alpha}(z, w)$, $\dfrac{\partial B_\delta}{\partial \bar{z}_\alpha}(z, w)$. Moreover, all the derivatives are bounded with respect to $w$ uniformly with $z$ running over a compact set. Also*

$$\frac{\partial B_\delta}{\partial \bar{w}_\gamma}(z, w) \to 0, \quad \frac{\partial B_\delta}{\partial \bar{z}_\alpha}(z, w) \to 0 \qquad as \quad |w'| \to \infty \quad or \quad |w_\delta| \to \infty.$$

*Actually the above holds for any kind of first order derivatives.*

**Proof Lemma 3.** The first two equations in (19) are obvious, so we have to show that we can enter with differentiation under the integral sign (the third equation). One way is to apply theorems about differentiation of integrals with parameter, for instance see S.M. Nikolsky [19], pp. 118 - 146, but the theorems from the book do not exactly apply to our case, and some small work is required. Another way is to prove that the differentiation can be moved inside the integral by using the Lebesgue Dominated Convergence Theorem. We choose the latter method.

We consider $\partial/\partial \bar{w}_\gamma$ as the Wirtinger derivative

$$\frac{\partial}{\partial \bar{w}_\gamma} = \frac{1}{2}\left(\frac{\partial}{\partial u_\gamma} + i\frac{\partial}{\partial v_\gamma}\right), \qquad \gamma = 1, \ldots, k.$$

Of course, it is enough to prove our request for the derivative $\partial/\partial u_\gamma$ since the calculations are the same for other derivatives: $\partial/\partial v_\gamma$, $\partial/\partial x_\alpha$, $\partial/\partial y_\alpha$.

In order to simplify our calculations, we use the notation

$$\begin{aligned}
[w]_{\gamma,\delta}(h, \zeta) &= (w_1, \ldots, w_{\gamma-1}, w_\gamma + h, w_{\gamma+1}, \ldots, w_{\delta-1}, w_\delta + \zeta, w_{\delta+1}, \ldots, w_k) &&\text{if } \gamma \neq \delta \\
[w]_\delta(h, \zeta) &= (w_1, \ldots, w_{\delta-1}, w_\delta + \zeta + h, w_{\delta+1}, \ldots, w_k) &&\text{if } \gamma = \delta
\end{aligned} \quad (20)$$

Since we will be making calculations for $\dfrac{\partial B_\delta}{\partial u_\gamma}$, therefore $h$ is real. First we assume that $\gamma \neq \delta$; we have the following sequence of equations:

$$\frac{B_\delta(z, [w]_\gamma(h, 0)) - B_\delta(z, w)}{h} =$$

$$= \frac{1}{h}\left[\frac{1}{2\pi i}\int_{\mathbb{C}} \frac{b_\delta(z, [w]_{\gamma,\delta}(h, \zeta))}{\zeta} d\zeta \wedge d\bar\zeta - \frac{1}{2\pi i}\int_{\mathbb{C}} \frac{b_\delta(z, [w]_{\gamma,\delta}(0, \zeta))}{\zeta} d\zeta \wedge d\bar\zeta\right]$$

$$= \frac{1}{2\pi i}\int_{\mathbb{C}} \frac{b_\delta(z, [w]_{\gamma,\delta}(h, \zeta)) - b_\delta(z, [w]_{\gamma,\delta}(0, \zeta))}{\zeta h} d\zeta \wedge d\bar\zeta \tag{21}$$

The functions $b_\delta(z, w)$, $\delta = 1, \ldots, k$, are complex-valued, we apply a modified Mean Value Theorem for vector-valued functions (see e.g., W. Rudin [21] p. 113, or S.M. Nikolsky [19] p. 119). If we consider the function

$$\mathbb{R} \ni h \longrightarrow b_\delta(z, [w]_{\gamma,\delta}(h, \zeta)) - b_\delta(z, [w]_{\gamma,\delta}(0, \zeta)) \in \mathbb{C}$$

for $|h|$ sufficiently small (without any loss of generality we can assume that $|h| \leq 1$), there exits $\widehat{h}$, $|\widehat{h}| < |h|$, such that

$$\left|b_\delta(z, [w]_{\gamma,\delta}(h, \zeta)) - b_\delta(z, [w]_{\gamma,\delta}(0, \zeta))\right| \leq |h|\left|\frac{\partial b_\delta}{\partial u_\gamma}(z, [w]_{\gamma,\delta}(\widehat{h}, \zeta))\right|$$

$$\leq |h| \sup_{-|h| \leq \tau \leq |h|} \left|\frac{\partial b_\delta}{\partial u_\gamma}(z, [w]_{\gamma,\delta}(\tau, \zeta))\right| \tag{22}$$

Taking into account the meaning of $[w]_{\gamma,\delta}(\tau, \zeta)$, and the assumption that the form $\omega(z, w)$ is decreasing along the fibers, we get the estimate

$$|h|\left|\frac{\partial b_\delta}{\partial u_\gamma}(z, [w]_{\gamma,\delta}(\tau, \zeta))\right| \leq$$

$$\leq \frac{C|h|}{1 + |w_1|^{1+\varepsilon} + \ldots + |w_\gamma + \tau|^{1+\varepsilon} + \ldots + |w_{\delta-1}|^{1+\varepsilon} + |w_\delta + \zeta|^{1+\varepsilon} + |w_{\delta+1}|^{1+\varepsilon} + \ldots + |w_k|^{1+\varepsilon}}$$

$$\leq \frac{C|h|}{1 + |w_\delta + \zeta|^{1+\varepsilon}}$$

$$\leq \frac{\widetilde{C}|h|}{1 + |\zeta|^{1+\varepsilon}} \quad \text{for possibly larger } \widetilde{C} \quad (\text{because } w_\delta \text{ is running over a compact set}). \tag{23}$$

Using the inequalities (22), (23), the integrand in (21) can be estimated

$$\left|\frac{b_\delta(z, [w]_{\gamma,\delta}(\zeta + h)) - b_\delta(z, [w]_{\gamma,\delta}(\zeta))}{\zeta h}\right| \leq \frac{\widetilde{C}|h|}{|h||\zeta|(1 + |\zeta|^{1+\varepsilon})} = \frac{\widetilde{C}}{|\zeta|(1 + |\zeta|^{1+\varepsilon})} \tag{24}$$

As we know from the end of Subsection 3.1, the function $1/[|\zeta|(1+|\zeta|^{1+\varepsilon})]$ is integrable over $\mathbb{C}$. From the Lebesgue Dominated Convergence Theorem, we can pass to the limit as $h \to 0$ under the integral sign in (21) and we get

$$\lim_{h\to 0}\frac{B_\delta(z,[w]_{\gamma,\delta}(h,0)) - B_\delta(z,w)}{h} = \frac{1}{2\pi i}\int_\mathbb{C}\left[\lim_{h\to 0}\frac{b_\delta(z,[w]_{\gamma,\delta}(h,\zeta)) - b_\delta(z,[w]_{\gamma,\delta}(0,\zeta))}{\zeta h}\right]d\zeta \wedge d\bar\zeta$$

or equivalently

$$\frac{\partial B_\delta}{\partial u_\gamma}(z,w) = \frac{1}{2\pi i}\int_\mathbb{C}\frac{1}{\zeta}\cdot\frac{\partial b_\delta}{\partial u_\gamma}(z,w_1,\ldots,w_{\delta-1},w_\delta+\zeta,w_{\delta+1},\ldots,w_k)\,d\zeta\wedge d\bar\zeta \qquad (25)$$

where $\gamma,\delta = 1,\ldots,k$, $\gamma\neq\delta$.

If $\gamma = \delta$, calculations are practically the same as in the case $\gamma \neq \delta$. The only difference is that instead of taking $w_\gamma + h$ and $w_\delta + \zeta$, we use $w_\delta + \zeta + h$. So the formula (25) also holds for $\gamma = \delta$.

Also the same way the proof of (25) goes if we differentiate with respect to $v_\gamma = \operatorname{Im} w_\gamma$ instead of $u_\gamma = \operatorname{Re} w_\gamma$. Because of that, we have

$$\frac{\partial B_\delta}{\partial \bar w_\gamma}(z,w) = \frac{1}{2\pi i}\int_\mathbb{C}\frac{1}{\zeta}\cdot\frac{\partial b_\delta}{\partial \bar w_\gamma}(z,w_1,\ldots,w_{\delta-1},w_\delta+\zeta,w_{\delta+1},\ldots,w_k)\,d\zeta\wedge d\bar\zeta \qquad (26)$$

In exactly the same way we prove the property, of entering with differentiation under the integral sign, for all other first-order derivatives. The last part of the lemma follows immediately from Section 3.2, where we replace $B_\delta(z,w)$ by $DB_\delta(z,w)$, where $D$ is any first-order derivative.

## 3.4 Using the compatibility conditions

As the next step, we use the compatibility conditions (4) in the integration below:

$$\begin{aligned}
\frac{\partial B_\delta}{\partial \bar w_\gamma}(z,w) &= \frac{1}{2\pi i}\int_\mathbb{C}\frac{1}{\zeta}\cdot\frac{\partial b_\delta}{\partial \bar w_\gamma}(z,w_1,\ldots,w_{\delta-1},w_\delta+\zeta,w_{\delta+1},\ldots,w_k)\,d\zeta\wedge d\bar\zeta \\
&= \frac{1}{2\pi i}\int_\mathbb{C}\frac{1}{\zeta}\cdot\frac{\partial b_\gamma}{\partial \bar w_\delta}(z,w_1,\ldots,w_{\delta-1},w_\delta+\zeta,w_{\delta+1},\ldots,w_k)\,d\zeta\wedge d\bar\zeta \\
&= \frac{1}{2\pi i}\int_\mathbb{C}\frac{1}{\zeta}\cdot\frac{\partial b_\gamma}{\partial \bar\zeta}(z,w_1,\ldots,w_{\delta-1},w_\delta+\zeta,w_{\delta+1},\ldots,w_k)\,d\zeta\wedge d\bar\zeta \qquad (27)
\end{aligned}$$

At this point we use the Bochner-Martinelli formula in the case of one variable. Here we recall this formula:

**Theorem (see [17], Corollary 1.1.5, p. 24).** *If $\Omega \subseteq \mathbb{C}$ is a bounded domain with $C^1$ boundary and if $f \in C^1(\overline{\Omega})$, then, for any $z_0 \in \Omega$,*

$$f(z_0) = \frac{1}{2\pi i} \int_{\partial \Omega} \frac{f(\zeta)}{\zeta - z_0} d\zeta + \frac{1}{2\pi i} \int_\Omega \frac{\partial f(\zeta)/\partial \overline{\zeta}}{\zeta - z_0} d\zeta \wedge d\overline{\zeta}.$$

In our case, we apply the above theorem for the domain

$$\Omega = D(0, R) = \{z \in \mathbb{C}; \ |z| < R\}, \quad point \quad z_0 = 0,$$

and the function

$$f(\zeta) = b_\gamma(z, w_1, \ldots, w_{\delta-1}, w_\delta + \zeta, w_{\delta+1}, \ldots, w_k).$$

The domain $\Omega$ and the function $f = f(\zeta)$ satisfy all the assumptions of the theorem, so we obtain

$$b_\gamma(z, w) = \frac{1}{2\pi i} \int_{\partial D(0,R)} \frac{b_\gamma(z, w_1, \ldots, w_{\delta-1}, w_\delta + \zeta, w_{\delta+1}, \ldots, w_k)}{\zeta} d\zeta +$$

$$+ \frac{1}{2\pi i} \int_{D(0,R)} \frac{1}{\zeta} \cdot \frac{\partial b_\gamma}{\partial \overline{\zeta}}(z, w_1, \ldots, w_{\delta-1}, w_\delta + \zeta, w_{\delta+1}, \ldots, w_k) \, d\zeta \wedge d\overline{\zeta} \quad (28)$$

By assumption, the functions $b_\gamma(z, w)$ are decreasing along the fibers, so we get the following estimates

$$\left| \frac{b_\gamma(z, w_1, \ldots, w_{\delta-1}, w_\delta + \zeta, w_{\delta+1}, \ldots, w_k)}{\zeta} \right| \leqslant \frac{C}{|\zeta| \cdot (1 + |w_\delta + \zeta|^{1+\varepsilon})} \leqslant \frac{C}{|\zeta| \cdot \left(1 + ||\zeta| - |w_\delta||^{1+\varepsilon}\right)} \quad (29)$$

$$\left| \frac{\frac{\partial b_\gamma}{\partial \overline{\zeta}}(z, w_1, \ldots, w_{\delta-1}, w_\delta + \zeta, w_{\delta+1}, \ldots, w_k)}{\zeta} \right| \leqslant \frac{C}{|\zeta| \cdot (1 + |w_\delta + \zeta|^{1+\varepsilon})} \leqslant \frac{C}{|\zeta| \cdot \left(1 + ||\zeta| - |w_\delta||^{1+\varepsilon}\right)} \quad (30)$$

The first integral in (28) can be estimated

$$\left| \frac{1}{2\pi i} \int_{\partial D(0,R)} \frac{b_\gamma(z, w_1, \ldots, w_{\delta-1}, w_\delta + \zeta, w_{\delta+1}, \ldots, w_k)}{\zeta} d\zeta \right| \leqslant$$

$$\leqslant \frac{1}{2\pi} \int_{\partial D(0,R)} \left( \frac{1}{R} \cdot \frac{C}{1 + |R - |w_\delta||^{1+\varepsilon}} \right) |d\zeta| = \frac{1}{2\pi} \int_0^{2\pi} \left( \frac{1}{R} \cdot \frac{C}{1 + |R - |w_\delta||^{1+\varepsilon}} \right) R \, d\theta =$$

$$= \frac{1}{2\pi} \int_0^{2\pi} \frac{C}{1 + |R - |w_\delta||^{1+\varepsilon}} d\theta = \frac{C}{1 + |R - |w_\delta||^{1+\varepsilon}} \longrightarrow 0 \quad as \quad R \to \infty \quad (31)$$

Next, we consider the second integral in (28). We note that the integral taken over the entire plane $\mathbb{C}$ exists because of the following estimates:

$$\left| \frac{1}{2\pi i} \int_{\mathbb{C}\setminus D(0,R)} \frac{1}{\zeta} \cdot \frac{\partial b_\gamma}{\partial \bar\zeta}(z, w_1, \ldots, w_{\delta-1}, \zeta + w_\delta, w_{\delta+1}, \ldots, w_k) \, d\zeta \wedge d\bar\zeta \right| \leq$$

$$\leq \frac{1}{2\pi} \int_{\mathbb{C}\setminus D(0,R)} \frac{1}{|\zeta|} \cdot \left| \frac{\partial b_\gamma}{\partial \bar\zeta}(z, w_1, \ldots, w_{\delta-1}, \zeta + w_\delta, w_{\delta+1}, \ldots, w_k) \right| |d\zeta \wedge d\bar\zeta| \leq$$

$$\leq \frac{1}{2\pi} \int_0^{2\pi} \int_R^\infty \frac{1}{\rho} \cdot \frac{C}{1 + \big|\rho - |w_\delta|\big|^{1+\varepsilon}} \cdot 2\rho \, d\rho \, d\theta \leq 2 \int_R^\infty \frac{C}{1 + \big|\rho - |w_\delta|\big|^{1+\varepsilon}} \, d\rho \longrightarrow 0 \quad \text{as} \quad R \to \infty \quad (32)$$

From the estimates (30) and (32) we see that the integral over the entire plane $\mathbb{C}$

$$\frac{1}{2\pi i} \int_\mathbb{C} \frac{1}{\zeta} \cdot \frac{\partial b_\gamma}{\partial \bar\zeta}(z, w_1, \ldots, w_{\delta-1}, \zeta + w_\delta, w_{\delta+1}, \ldots, w_k) \, d\zeta \wedge d\bar\zeta$$

is convergent.

Finally, combining (28) and the above estimates (30) and (32), and taking the limit when $R \to \infty$, we get the formula

$$b_\gamma(z, w) = \frac{1}{2\pi i} \int_\mathbb{C} \frac{1}{\zeta} \cdot \frac{\partial b_\gamma}{\partial \bar\zeta}(z, w_1, \ldots, w_{\delta-1}, w_\delta + \zeta, w_{\delta+1}, \ldots, w_k) \, d\zeta \wedge d\bar\zeta \quad (33)$$

Comparing the formulas (27) and (33), we obtain

$$\frac{\partial B_\delta}{\partial \bar w_\gamma}(z, w) = b_\gamma(z, w), \qquad \gamma, \delta = 1, \ldots, k, \quad (34)$$

or equivalently

$$\bar\partial_w B_\delta(z, w) = b(z, w).$$

We note that the right-hand side of the equation does not depend on $\delta$. Take two functions $B_\delta(z, w)$ and $B_{\delta'}(z, w)$ which satisfy

$$\bar\partial_w B_\delta(z, w) = \bar\partial_w B_{\delta'}(z, w) \quad \text{or} \quad \bar\partial_w \left[ B_\delta(z, w) - B_{\delta'}(z, w) \right] = 0$$

which implies

$$w \longrightarrow B_\delta(z, w) - B_{\delta'}(z, w) \quad \text{is holomorphic with respect to } w.$$

Since $B_\delta(z, w)$ and $B_{\delta'}(z, w)$ are bounded with respect to $w$ and

$$\lim_{\|w'\| \to \infty} B_\delta(z, w) = 0, \quad \lim_{|w_\delta| \to \infty} B_\delta(z, w) = 0, \quad \lim_{\|w'\| \to \infty} B_{\delta'}(z, w) = 0, \quad \lim_{|w_\delta| \to \infty} B_{\delta'}(z, w) = 0$$

we obtain, by the Liouville's Theorem, that

$$B_\delta(z, w) = B_{\delta'}(z, w) \qquad \text{for any} \quad \delta, \delta' = 1, \ldots, k. \quad (35)$$

Therefore dropping the index $\delta$ is justified,

$$B(z,w) = \frac{1}{2\pi i}\int_{\mathbb{C}} \frac{b_\delta(z,w_1,\ldots,w_{\delta-1},\zeta,w_{\delta+1},\ldots,w_k)}{\zeta - w_\delta}\, d\zeta \wedge d\bar\zeta$$

$$= \frac{1}{2\pi i}\int_{\mathbb{C}} \frac{b_\delta(z,w_1,\ldots,w_{\delta-1},w_\delta+\zeta,w_{\delta+1},\ldots,w_k)}{\zeta}\, d\zeta \wedge d\bar\zeta$$

## 3.5 Globally defined form on the fiber bundle

In this subsection we repeat the argument given in [7], pp. 562 - 563, that the functions $B(z,w) = B_s(z,w)$, $s \in S$, defined in each coordinate domain $X_s \times \mathbb{C}^k$, actually can be put together and give a global function.

In another chart, say $X_t \simeq U' \times \mathbb{C}^k$, the form $\omega'(z',w')$ coincides with $\omega(z,w)$ on the intersection $(U \cap U') \times \mathbb{C}^k$ and we have the corresponding functions

$$B'(z',w') = \frac{1}{2\pi i}\int_{\mathbb{C}} \frac{b'_\delta(z',w'_1,\ldots,w'_{\delta-1},\zeta,w'_{\delta+1},\ldots,w'_k)}{\zeta - w'_\delta}\, d\zeta \wedge d\bar\zeta$$

Of course, also we have

$$B' = B'_\gamma = B'_\delta \quad \text{and} \quad \frac{\partial}{\partial \overline{w'_\gamma}} B'(z',w') = b'_\gamma(z',w'), \quad \gamma,\delta = 1,\ldots,k, \tag{36}$$

We will prove that $B(z,w) = B'(z',w')$ on the domain where they both are defined. To avoid a confusion, here $w'$ is just another point $w$, i.e., $w' = (w'_1,\ldots,w'_k)$.

We have

$$d\bar z'_\alpha = \sum_{\beta=1}^n \frac{\overline{\partial f_\alpha}}{\partial \bar z_\beta}(z,w)\, d\bar z_\beta \qquad \alpha = 1,\ldots,n$$

$$d\overline{w'_\gamma} = \sum_{\beta=1}^n \frac{\partial \overline{g_\gamma}}{\partial \bar z_\beta}(z,w)\, d\bar z_\beta + \sum_{\delta=1}^k \frac{\partial \overline{g_\gamma}}{\partial \overline{w_\delta}}(z,w)\, d\overline{w_\delta} \qquad \gamma = 1,\ldots,k.$$

Now we calculate the change of coefficients of $d\overline{w_\gamma}$ in the form $\omega$ from (2). We have

$$\sum_{\gamma=1}^k b'_\gamma(z',w')\, d\overline{w'_\gamma} = \sum_{\gamma=1}^k \left[ b'_\gamma(z',w') \left( \sum_{\beta=1}^n \frac{\partial \overline{g_\gamma}}{\partial \bar z_\beta}(z,w)\, d\bar z_\beta + \sum_{\delta=1}^k \frac{\partial \overline{g_\gamma}}{\partial \overline{w_\delta}}(z,w)\, d\overline{w_\delta} \right) \right]$$

$$= \sum_{\beta=1}^n \left[ \sum_{\gamma=1}^k b'_\gamma(z',w') \frac{\partial \overline{g_\gamma}}{\partial \bar z_\beta}(z,w) \right] d\bar z_\beta + \sum_{\delta=1}^k \underbrace{\left[ \sum_{\gamma=1}^k b'_\gamma(z',w') \frac{\partial \overline{g_\gamma}}{\partial \overline{w_\delta}}(z,w) \right]}_{b_\delta(z,w)} d\overline{w_\delta}.$$

So we got

$$b_\delta(z,w) = \sum_{\gamma=1}^{k} b'_\gamma(f(z), g(z,w)) \frac{\partial \overline{g}_\gamma}{\partial \overline{w}_\delta}(z,w) \tag{37}$$

To simplify calculations, we use the notation

$$[w]_\delta(\zeta) = (w_1, \ldots, w_{\delta-1}, \zeta, w_{\delta+1}, \ldots, w_k).$$

Using (36) we have

$$\frac{\partial}{\partial \overline{\zeta}}\left[B'(f(z), g(z,[w]_\delta(\zeta)))\right] = \sum_{\gamma=1}^{k} b'_\gamma(f(z), g(z,[w]_\delta(\zeta))) \frac{\partial \overline{g}_\gamma}{\partial \overline{w}_\delta}(z, [w]_\delta(\zeta)) \tag{38}$$

Then we calculate the integral (6) - (7), using (37) in the second equation below and (38) in the third one, we obtain

$$\begin{aligned}
B(z,w) &= \frac{1}{2\pi i} \int_{\mathbb{C}} \frac{b_\delta(z, [w]_\delta(\zeta))}{\zeta - w_\delta} \, d\zeta \wedge d\overline{\zeta} \\
&= \frac{1}{2\pi i} \int_{\mathbb{C}} \frac{1}{\zeta - w_\delta} \sum_{\gamma=1}^{k} b'_\gamma(f(z), g(z,[w]_\delta(\zeta))) \frac{\partial \overline{g}_\gamma}{\partial \overline{w}_\delta}(z, [w]_\delta(\zeta)) \, d\zeta \wedge d\overline{\zeta} \\
&= \frac{1}{2\pi i} \int_{\mathbb{C}} \frac{1}{\zeta - w_\delta} \cdot \frac{\partial}{\partial \overline{\zeta}}\left[B'(f(z), g(z,[w]_\delta(\zeta)))\right] d\zeta \wedge d\overline{\zeta} \\
&= \frac{1}{2\pi i} \int_{\mathbb{C}} \frac{1}{\zeta - w_\delta} \cdot \frac{\partial}{\partial \overline{\zeta}}\left[B'\left(f(z), g(z, w_1, \ldots, w_{\delta-1}, \zeta, w_{\delta+1}, \ldots, w_k)\right)\right] d\zeta \wedge d\overline{\zeta} \\
&= B'(f(z), g(z,w))
\end{aligned}$$

what we wanted to prove.

## 3.6 End of the proof of the theorem

We recall the main points what we proved: In each coordinate map $X_s \sim U_s \times \mathbb{C}^k$, we have a solution $B(z,w) = B_s(z,w)$ that satisfies

1. $B(z,w)$ is bounded with respect to $w$ and $B(z,w) \longrightarrow 0$ as $\|w'\| \to \infty$ or $|w_\delta| \to \infty$ for some $\delta$, uniformly in each chart $X_s = U_s \times \mathbb{C}^k$ with respect to $z$ running over a compact set.

2. $\dfrac{\partial B}{\partial \bar{w}_\gamma}(z,w) = b_\gamma(z,w), \quad \gamma = 1,\ldots,k, \quad \text{or equivalently} \quad \bar{\partial}_w B(z,w) = b(z,w).$

3. $\dfrac{\partial B}{\partial \bar{z}_\alpha}(z,w)$ *exists for any* $\alpha = 1,\ldots,n,$ *and* $\dfrac{\partial B}{\partial \bar{z}_\alpha}(z,w) \longrightarrow 0$ *as* $\|w'\| \to \infty$ *or* $|w_\delta| \to \infty$
*for* $\delta = 1,\ldots,k,$ *uniformly in each chart* $X_s = U_s \times \mathbb{C}^k$ *with respect to* $z$ *running over a compact set.*

We write $\bar{\partial}_z B$ in terms of a combination of $d\bar{z}_\alpha$,

$$\bar{\partial}_z B(z,w) = \hat{a}_1(z,w)\, d\bar{z}_1 + \ldots + \hat{a}_n(z,w)\, d\bar{z}_n, \qquad \hat{a}_\alpha(z,w) = \dfrac{\partial B}{\partial \bar{z}_\alpha}(z,w).$$

Clearly, the both $(0,1)$-forms, $\bar{\partial} B$ and $\omega$,

$$\bar{\partial} B(z,w) = \bar{\partial}_z B(z,w) + \bar{\partial}_w B(z,w) \qquad \text{and} \qquad \omega(z,w) = a(z,w) + b(z,w)$$

are $\bar{\partial}$-closed and

$$\bar{\partial}_w B(z,w) = b(z,w)$$

Consequently, using a particular case of the compatibility equations (4), namely

$$\dfrac{\partial \hat{a}_\alpha}{\partial \bar{w}_\gamma}(z,w) = \dfrac{\partial b_\gamma}{\partial \bar{z}_\alpha}(z,w) \qquad \text{and} \qquad \dfrac{\partial a_\alpha}{\partial \bar{w}_\gamma}(z,w) = \dfrac{\partial b_\gamma}{\partial \bar{z}_\alpha}(z,w), \quad \alpha = 1,\ldots,n,\ \gamma = 1,\ldots,k,$$

and from here

$$\dfrac{\partial \hat{a}_\alpha}{\partial \bar{w}_\gamma}(z,w) = \dfrac{\partial a_\alpha}{\partial \bar{w}_\gamma}(z,w), \qquad \alpha = 1,\ldots,n,\ \gamma = 1,\ldots,k$$

or

$$\dfrac{\partial(\hat{a}_\alpha - a_\alpha)}{\partial \bar{w}_\gamma}(z,w) = 0 \qquad \text{for} \quad \alpha = 1,\ldots,n,\ \gamma = 1,\ldots,k$$

therefore
$$\hat{a}_\alpha(z,w) - a_\alpha(z,w), \quad \alpha = 1,\ldots,n, \quad \text{are holomorphic with respect to } w.$$

Since the both functions $a_\alpha(z,w)$ and $\hat{a}_\alpha(z,w)$ are bounded with respect to $w$ and

$$a_\alpha(z,w) \longrightarrow 0, \qquad \hat{a}_\alpha(z,w) \longrightarrow 0 \qquad \text{as} \quad \|w'\| \to \infty, \text{ or } |w_\delta| \to \infty \text{ for } \delta = 1,\ldots,k,$$

from the Liouville's Theorem,
$$\hat{a}_\alpha(z,w) - a_\alpha(z,w) = \text{const}_\alpha = 0.$$

We have
$$\hat{a}_\alpha(z,w) = a_\alpha(z,w), \qquad \alpha = 1,\ldots,n, \qquad \text{or equivalently} \quad \bar{\partial}_z B(z,w) = a(z,w)$$

which together with $\bar{\partial}_w B(z,w) = b(z,w)$, yields

$$\bar{\partial} B(z,w) = \omega(z,w).$$

The main theorem is proved.

# 4 Application of the Main Theorem to bundles over $\mathbb{CP}^1$

In this section we show a simple application of the main theorem in the situation that other results cannot be applied. Let $\mathbb{CP}^2$ be the 2-dimensional complex projective space. We know that $\mathbb{C}^2$ can be embedded into $\mathbb{CP}^2$ and we can write
$$\mathbb{CP}^2 = \mathbb{C}^2 \cup \mathbb{CP}^1_\infty$$
i.e., we can compactify $\mathbb{C}^2$ by adding $\mathbb{CP}^1$ at infinity (see [7] or [11]).

Similarly, if we remove a point $p$ from $\mathbb{CP}^2$, then we obtain a vector bundle $E$ over $\mathbb{CP}^1$ with one-complex dimensional fiber (see [7]).

We will make some comments about the latter case because it is more interesting and directly the main theorem can be applied. Namely, let $\omega$ be a closed $(0,1)$ form defined on $\mathbb{CP}^2$ that vanishes to order two or more at the point $p$, mentioned above. This form $\omega$ can be considered as a form defined on the vector bundle $E$. Of course, this form is closed and, applying the main theorem, there is a solution of the equation $\bar\partial \sigma = \omega$. Moreover the solution $\sigma$ vanishes along the fibers of $E$, however, not necessarily has compact support along fibers. From the main theorem also we can deduce how fast the function $\omega$ is approaching zero if the argument goes to "infinity" along the fibers. Of course, this speed depends on the order of decreasing of the form $\omega$.

Małgorzata Urlińska
Warsaw University of Technology
00-662 Warsaw, Poland
E-mail: m.urlinska@gmail.com